\documentclass{article}
\usepackage{graphicx} % Required for inserting images
\usepackage{amsmath, amssymb}
\usepackage[portrait, margin=3cm]{geometry}

\usepackage[square,numbers]{natbib}
\newcommand{\ba}{\ensuremath{\mathbf{a}}}
\newcommand{\bb}{\ensuremath{\mathbf{b}}}

\newcommand{\bh}{\ensuremath{\mathbf{h}}}
\newcommand{\bi}{\ensuremath{\mathbf{i}}}
\newcommand{\bj}{\ensuremath{\mathbf{j}}}
\newcommand{\bk}{\ensuremath{\mathbf{k}}}

\newcommand{\bx}{\ensuremath{\mathbf{x}}}

\newcommand{\bA}{\ensuremath{\mathbf{A}}}
\newcommand{\bB}{\ensuremath{\mathbf{B}}}

\newcommand{\bX}{\ensuremath{\mathbf{X}}}

\newcommand{\etal}{{\it{et al}}.}

\newcommand{\Zint}{\mathbb{Z}}

\usepackage{bm}

\title{A Complete Graphic Statics for Rigid-Jointed 3D Frames. \\Part 2: Homology of loops}
\author{Allan McRobie, fam20@cam.ac.uk \\
Dept. of Engineering, Cambridge University, CB2 1PZ\\
https://orcid.org/0000-0002-6610-5927}
\begin{document}

\maketitle

\section*{\begin{center} Abstract \end{center}}
This paper, Part 2 of a sequence, seeks to extend the reach of graphic statics by describing the forces and moments in any 3D rigid-jointed frame structure in terms of cell complexes using homology theory of algebraic topology. Graphic statics provides a highly geometric way to represent the equilibrium in bar structures. Unlike traditional matrix-based linear structural analysis which represents a structure as a set of nodes connected by bars,  graphic statics imagines that the bar network defines a variety of 
higher-dimensional objects (polygonal faces, polyhedral cells, polytopes,~$\ldots$). These objects are related to piecewise-linear stress functions, the liftings of Maxwell, Rankine or Cremona. The requirement for such stress-functions to be plane-faced places a major limitation on the set of structures that can be analysed, as in many structures the spaces between bars do not correspond to flat polygonal regions.

The CW-complexes of cellular homology provide a far-reaching generalisation of geometric notions such as polygons, polyhedra and polytopes, and their use here
removes the requirement that spaces between bars must be flat. 
 Here we demonstrate how any frame structure with bar-like members can be decomposed into a union of closed loops, each consisting of a closed circuit of bars.  For general structures these loops are general closed space curves which cannot be spanned by flat polygons.  Using chains of CW-complexes makes the new theory 
applicable to a much richer set of structural geometries. 

Unlike most descriptions of graphic statics, this approach is not restricted to purely axial forces. Shear forces, bending moments and torsional moments are included naturally, as described in Part 1 of this sequence of papers. Later papers in the sequence will extend the approach to displacements, rotations and Virtual Work, and will give greater detail on how the loop formalism may be ``lifted'' to involve higher dimensional CW-complexes.

\section{Introduction}
Graphic statics~\cite{allen2010, MaxwellReader, maxwell1864a, maxwell1870a} has its limitations and the main purpose of this paper is to generalise the approach such that a wider set of structural problems can be addressed. This is achieved primarily by replacing the familiar geometrical concepts of polygon, polyhedron and polytope with the more general notion of CW-complexes, as used in algebraic topology. The key objects in the new description are loops in a 4D space. Sets of loops are used to represent the structural frame and a set of dual loops will represent the forces and moments within the members of the frame. The ability to represent moments is a key strength of this new approach. Later papers will show how other loops can represent kinematic variables such as displacements and rotations.

Maxwell's classic 1864 paper~\cite{maxwell1864a} on reciprocal diagrams makes the telling observation that if the form diagram of a 2D truss is the 2D projection of a plane-faced polyhedron then the truss can carry a state of axial self-stress. This connection is consolidated in his much longer 1870 paper~\cite{maxwell1870a}, where he identifies the polyhedron as a discrete Airy stress function, with all curvature (and hence all stresses) concentrated along the sharp edges of the polyhedron. Maxwell proposed various stress functions to deal with 3D stresses. Here we focus on the one implicit in his ``Diagram of Stress'' construction which is a Legendre transform between body space and stress space~\cite{McRobieArXiv1, LegendreReader, McRobieZurich}, applicable in any dimension. 
Part 1 of this sequence of papers~\cite{McRobieArXiv1} showed how the forces and moments within a simple structure consisting of a single closed loop of bars could be represented by a dual loop.
The dual loop exists in 4D space - this being the usual 3D space with an extra dimension for the stress function. The areas of three of its six projections give the axial and two shear forces whilst the other three projections give the components of the total moment about the origin. More general structures can then be represented as assemblages of such loops.

Using elementary homology theory, any frame structure can be decomposed into a union of loops, and the forces and moments can be represented by a union of dual loops. The resulting description, presented in Section 2 here, is complete: it can represent any state of self-stress in any frame structure. For the theory to align more closely with graphic statics, later papers in this series include ``liftings''. These are the higher dimensional objects (such as the faces of polyhedra and the cells of polytopes) which are associated with the stress functions of Maxwell's Diagram of Stress. It is these higher dimensional objects that distinguish graphic statics from the traditional linear analysis description of structural frame analysis. Note that the only structural objects in the Maxwell-Calladine count~\cite{Calladine1978} of the self-stresses and mechanisms of a structural truss are the nodes and bars. There is no mention of faces. This is in line with standard linear matrix structural analysis.  By contrast, faces are included in the count for the Euler characteristic of a polyhedral Airy stress function which plays a central role in graphic statics. Note that the faces do not exist. They are ideations - the analyst sees patterns -  triangles, quads and the like - in the gaps between the bars, and then goes on to make demands such as requiring these empty spaces to remain plane. It is this focus on the geometry and topology of the spaces between the bars that is so specific to graphic statics. The benefit of this perspective is that there is great explanatory power to be gained from seeing a truss as the projection of a polyhedron. Given a complicated 2D truss, a singular-value decomposition of its equilibrium matrix can reveal the possible states of self-stress and the possible inextensional mechanisms. These same features can often also be extracted by visual inspection of the polyhedral stress function, without recourse to matrix manipulation. This is the beauty of graphic statics: the analyst can often see at a glance how the structural forces are organised. In his 1864 paper~\cite{maxwell1864a}, Maxwell was cautious about extending this argument to 3D structures, stating that 
``the mechanical interest of reciprocal figures
in space rapidly diminishes with their complexity''. Whilst this has often been used as a rationale for not pursuing 3D reciprocal figures at all, here we use it a guide: 3D reciprocals will be investigated - and with moments - but we shall always aim to minimise the complexity whilst doing so.

\subsection{Other work}
Many results in this paper have been published previously by the author. There have been presentations at conferences of the International Association of Shell and Spatial Structures (IASS), particularly in Amsterdam 2015 ~\cite{McRobie2015}, Boston 2018~\cite{McRobieBoston}, Barcelona 2019~\cite{McRobieBarcelona} and Zurich 2024~\cite{McRobieZurich}. Part 1 of this sequence of papers provides further details of associated work by Karpenkov~\etal~\cite{Karpenkov2023} and Baranyai~\cite{Baranyai2024}. There is also work by Cooperband and co-authors using sheaf cohomology~\cite{Zoe3, Zoe2, ZoePhD} to study self-stresses and mechanisms in trusses and moment-resisting frames. Although there are significant overlaps in the mathematics - homology and cohomology theory - one difference is that - except for a few 2D examples -  the co-sheaf formalism involves no ``lifting''. There is no consideration (or imagining of)  higher dimensional objects like faces, cells or polytopes. Instead, the co-sheaf formalism is used to ascribe moments and forces to nodes and bars. One result~\cite{Zoe3} is a generalisation of the Maxwell-Calladine count for pin-jointed trusses to the case of moment-resisting frames.

\section{Homology of loops}
We begin with some preliminaries of homology theory. Together with homotopy theory, homology theory provides ways for describing the topological properties of objects. One advantage of homology theory is that it is a commutative theory, and many of its constituent parts are amenable to analysis using linear algebra. This property will provide a link between the highly geometrical constructions of graphic statics and the more familiar matrix-based methods of linear structural analysis employed by most practising structural engineers. There are at least three approaches to homology in  algebraic topology: simplicial, singular and cellular. The first two are based on simplices, (i.e.~triangles, tetrahedra and their higher dimensional generalisations). They require cells to be attached to each other according to strict rules, with vertices, edges, faces, etc.~numbered in specific orders. Those rules are helpful when constructing objects but are not sufficiently general for the purposes here, when the object - the structure - is already given. We shall thus adopt cellular homology, which uses CW-complexes instead of simplices.

Homology theory belongs to algebraic topology where it is typically used in the topological classification of objects. The aims here diverge from this, in that we are given a 3D frame structure and we seek to construct valid dual objects which represent the possible stress resultants within it. Despite the different objectives, the language and techniques of homology theory provide the very tools needed to generalise graphic statics.  We thus first describe the rather abstract mathematical machinery which, although elementary for topologists, may be unfamiliar to most structural engineers.

\begin{figure*}
\begin{center}
\includegraphics[width=\textwidth]{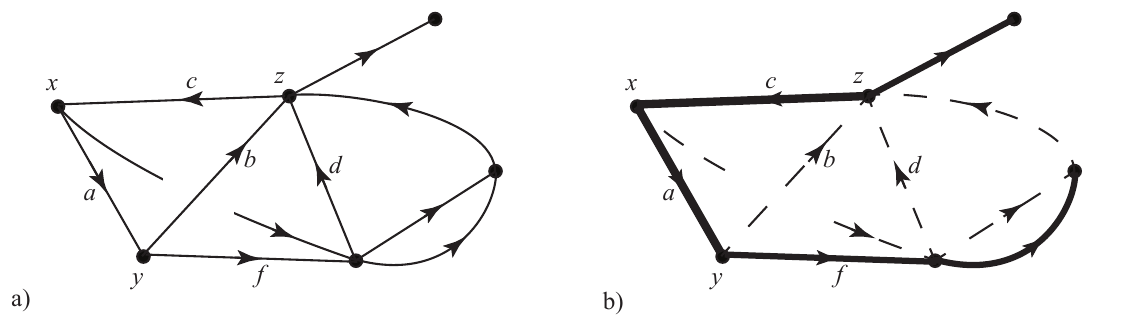}
\end{center}
\caption{a) A structure in 3D labelled as a graph, with nodes $x,y,z,\ldots$ and directed edges $a,b,c, \ldots$. b) A spanning tree for the graph.
}
\label{Fig1Graph}
\end{figure*}

The first step is to interpret a 3D rigid-jointed frame as a directed graph $X$. That is, $X$ consists of $v$ nodes connected by $e$ edges (the structural bars). The orientation of any bar (its direction in the directed graph) may be chosen freely. In the example of Fig.~\ref{Fig1Graph}$a$, there are six nodes, $\lbrace x, y, z, \ldots \rbrace$ and ten bars $\lbrace a, b, c, \ldots \rbrace$ with directions chosen as shown.

Next, the set of nodes is equipped with a group structure. That is, we define the free Abelian group $C_0$ generated by the vertices $\lbrace x, y, z, \ldots \rbrace$. Any integer linear combination of vertices, an object of the form $\alpha x + \beta y + \gamma z + \ldots$ (where $\alpha, \beta, \gamma \ldots \in \Zint$, the integers), is an element of this group. The group operation is given by ordinary addition of the integer coefficients. As example, the object $3x+4y-2z$ is an element of the group. That the group is Abelian means that it is commutative, such that for example $x+y = y+x$, meaning that the order in which the vertices are listed does not matter. The identity of the group is the zero element, written $0$, which is the object where all coefficients  $\alpha, \beta, \gamma \ldots $ equal zero. Note that the labels $x,y,z,\ldots$ denote the nodes themselves, and not their coordinates. Whilst the creation of such a group may be commonplace in topology, it is unusual in structural engineering. One rarely hears phrases such as ``take three of these nodes and subtract five of those nodes''. Nevertheless this ability to add and subtract objects in this manner will be of great use, especially as we progress to describe oriented bars.

For the bars, we define the free Abelian group $C_1$ generated by the directed edges $\lbrace a, b, c, \ldots \rbrace$. Any integer linear combination of directed edges, such as $5a+ 3b-7c$, is an element of this group. The group operation and the identity element are similar to $C_0$. Importantly, the Abelian nature of the group means that the order in which paths are listed does not matter. That is, the object $a+b-d$ is the same as the object $a-d+b$. There is no need for $b$ to follow $a$.  This is unlike homotopy theory, where paths are concatenated sequentially, and the order matters. Note that $a - d$ is also a valid element, even though its two component parts do not connect. 

We now define a map $\partial$ which takes elements of the edge group $C_1$ to elements of the node group $C_0$. This is a group homomorphism. It is the {\it boundary operator}, and it is defined by its action on the generators $\lbrace a, b, c, \ldots \rbrace$ of $C_1$.
Specifically, the action of $\partial$ on any directed edge (such as $a$) is given by its final vertex minus its initial vertex.
In this example,
\begin{eqnarray}
\partial: C_1 & \rightarrow & C_0 \\
\partial(a) & = & y - x
\end{eqnarray}

A {\it cycle} in $C_1$ is defined to be any element $t \in  C_1$ such that $\partial(t) = 0$, where $0$ is the identity element of $C_0$. A cycle then, is any set of directed edges (bars) which has no boundary. Such cycles are the loops that will be the primary objects of our graphical construction. In the example, it is clear that the object $t= a+b+c$ is a cycle. Perhaps less obvious would be that an element such as $u = a+2b+c-d-f$ is also a cycle. That the boundary $\partial(u) = 0$ may be checked by adding start and end nodes. Whilst we may be inclined to think of a loop (a cycle) as a closed path obtained by going for a walk around a set of edges in order (such as starting at vertex x, and following the path $u = a+b-d-f+b +c$) we may equally, since the group is Abelian,  define this loop by listing the relevant directed edges in any order.

The cycles are thus the kernel of the map $\partial:C_1 \rightarrow C_0$. That is, they are the elements of the edge group $C_1$ which map to the zero element of the node group, $0 \in C_0$. For simplicity, we shall assume there is only one structure. To make progress, choose a set of basis cycles which span this kernel. The most convenient way to achieve this is via a {\it spanning tree}. By definition, a tree has no loops, and we select any connected tree which is a subgraph of our graph $X$ and which visits every vertex.  An example is given in Fig.~\ref{Fig1Graph}$b$. By definition, a tree cannot contain any cycle, and it follows readily that every edge that is not in the tree defines an independent cycle. (Proof: any non-tree edge has a start and end node, every node is in the tree and there is a path from the end node back to the start node that passes through the tree, and which therefore forms a cycle when combined with the original edge.) In the example of Fig.~\ref{Fig1Graph}$b$, bar $d$ is not in the tree, and it defines the cycle $a+c+d+f$, which is a basis element for the kernel of $\partial: C_1 \rightarrow C_0$.

Trivially, a spanning tree has $v-1$ edges, and the number of edges not in the tree is thus $e-(v-1) = e-v+1$. There are thus $e-v+1$ independent cycles required to form the basis of the group of all cycles, which is the kernel of $\partial:C_1 \rightarrow C_0$. In the example of Fig.~\ref{Fig1Graph}$b$, there are five basis cycles, one for each bar not in the spanning tree.

\subsection{Dual loops}

\begin{figure*}[ht]
\begin{center}
\includegraphics[width=0.8\textwidth]{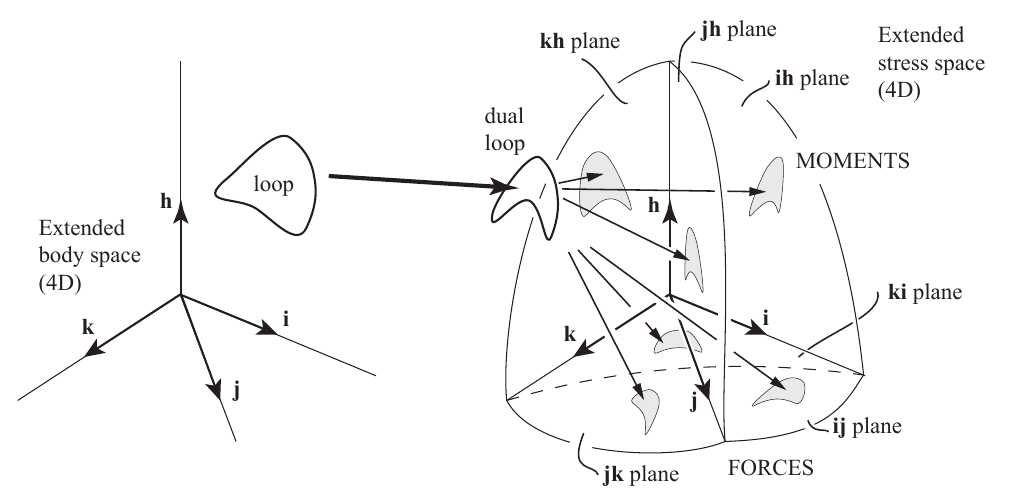}
\end{center}
\caption{A structural loop and its dual loop. The extended body and stress spaces have four dimensions, the usual three ($\bi, \bj$ and $\bk$) plus an extra one $\bh$ for the stress function. There are six basis bivector planes in the 4D stress space. Projections of the dual loop onto the planes $\bi\bj$, $\bj\bk$ and $\bk\bi$ give the force components.  Projections onto the planes $\bi\bh$, $\bj\bh$ and $\bk\bh$ give the moment components. 
}
\label{SixProjections}
\end{figure*} 
\begin{figure*}[ht]
\begin{center}
\includegraphics[width=\textwidth]{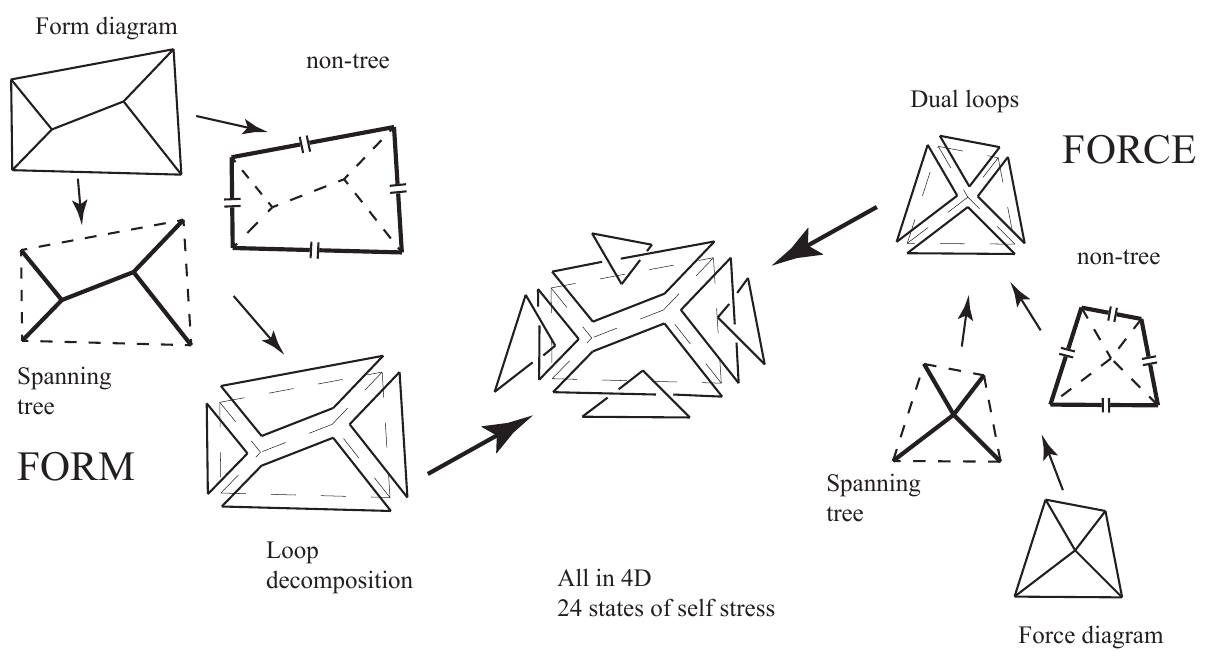}
\end{center}
\caption{Schematic illustrating a 3D frame structure with six nodes and nine bars being decomposed into four basis loops. The force diagram thus consists of four dual loops. }
\label{LoopDecomposition}
\end{figure*}

As described in Part 1 of this sequence of papers~\cite{McRobieArXiv1}, the forces and moments within a structural loop may be defined via dual loops in the 4D extended stress space. This follows from the Legendre transform underlying Maxwell's Diagram of Stress~\cite{maxwell1870a}. The oriented (bivector)  area of any dual loop in the 4D extended stress space has six components (see Fig.~\ref{SixProjections}). Three of these, on the usual $\bi \wedge \bj$, $\bj \wedge \bk$ and $\bk \wedge \bi$ planes of the usual 3D space, can be used to define the three components of force. There are three other projected areas, on the planes $\bi \wedge \bh$, $\bj \wedge \bh$ and $\bk \wedge \bh$. These planes
involve $\bh$, the direction in 4D which is associated with the stress function. These additional three projections provide the three components of the total moment within each basis structural loop.

Recall that the total moment includes both the internal moments (with a torsional and two bending moment components) and the moments due to the forces at the cut acting at a lever arm about the origin~\cite{McRobieArXiv1}. For a structure consisting of a single loop, the force and the total moment is thus constant around the loop, there being no applied forces. A single dual loop can thus represent the total stress resultant of force and total moment at any point on the structural loop. This is unlike the Corsican sum construction~\cite{McRobieRSOS1} where the dual object evolves around the loop. As shown in~\cite{McRobieArXiv1}, that construction corresponds to local bending and torsional moments, which vary from point to point around the loop.

Already then, this loop formalism can describe any state of self-stress in any 3D frame structure (and by inclusion, in any 2D frame also). First a spanning tree is identified to define a set of basis loops, and dual to each such loop we may define a dual loop in the 4D stress space to define the six components of the total stress resultant at a cut on each structural basis loop. The state of self-stress is then fully defined. Moreover, this description can define all possible states of self-stress in the frame.

Fig.~\ref{LoopDecomposition} provides a schematic to illustrate. A 3D frame structure with six nodes and nine bars is decomposed into four loops. The force diagram therefore consists of four dual loops. Any four dual loops arranged in any manner will define a state of self-stress. The stress resultant in each non-tree bar is given by the projected areas of the single loop that is dual to the basis structural loop defined by the non-tree bar. The stress resultant within a structural bar within the spanning tree is given by the projected areas of the union of those loops dual to the structural loops of which the bar is an edge. 

%In Fig.~\ref{LoopDecomposition}, the force loops have been assembled to create a force diagram a spanning tree has been drawn. This is not strictly necessary: the force diagram has been created as a set of loops, thus there is no need to decompose it again.

Any four force loops in any configuration would suffice to define a state of self-stress. In Fig.~\ref{LoopDecomposition}, force loops have been chosen that assemble into a comparatively neat  unified force diagram, but such neatness is unnecessary.

\subsection{Structural interpretation}
The description already contains structural engineering insight. If every bar not in the spanning tree is cut somewhere along its length, the frame becomes a 
statically-determinate tree structure containing no cycles. We may apply six independent components of the stress resultant at each cut, equal and opposite on the opposing cut faces. Each of these self-equilibrating release forces generates a self-stress in the loop (the cycle) defined by the cut bar. There are thus $6(e-v+1)$ independent states of self-stress in a 3D frame with $v$ nodes and $e$ bars. This is a well-known result in structural engineering. The Force Method for analysing statically indeterminate structures follows just this rationale: as is taught at engineering undergraduate level, a frame structure is made statically determinate (i.e.~a tree) by making the necessary number of cuts, and release forces and moments are then applied at each cut. Any state of stress in the frame is then the sum of the unique internal stresses in equilibrium with the external loads applied to the statically-determinate tree structure, plus any combination of the unique internal stress distributions which are in equilibrium with each of the release forces applied to the tree.

Traditionally, the Force Method continues by assuming linear elastic material behaviour and solving for the release forces and moments necessary to close the cuts. That is a compatibility requirement, and it requires knowledge of the material properties. Such considerations of compatibility and material law do not concern us at this stage. Here the focus is on equilibrium, and as described in Part 1, there are no external loads, only self-stresses. (External loads can always be considered to be applied by bars from an adjacent loading frame. One then studies states of self-stress in the combined system of structure and loading frame.)

%Note that the above counting rule applies to a single, simply-connected structure, but readily generalises to  $6(e-v+N)$ for $N$ separate structures. Whilst this may appear to be of little benefit, as structures are usually analysed one at a time, the ability to readily generalise the methods to multiple disjoint objects may prove to be of value when constructing the dual objects, the force diagrams. In such cases, disjoint objects may be the norm.

There is considerable freedom: any spanning tree may be chosen, but more importantly, dual to any structural loop, any choice of dual loop is admissible. There are numerous benefits to such freedom. 
\begin{itemize}
\item The description is complete. Any state of self-stress in any 2D or 3D frame structure can be represented.
\item The inclusion of moments in the description greatly expands the remit of graphic statics, which has traditionally been limited to truss structures capable of carrying axial forces only.
\item The polyhedral and polytopic stress functions in Rankine 3D reciprocals are subsumed within this description, but they have been considerably generalised here. For example, there is no need for loops to be the boundaries of plane polygonal faces. 
\end{itemize}

It is the representation of the internal forces and moments by loops that differentiates the approach in this paper from the methods of traditional structural analysis. Making a sufficient number of releases to decompose a frame structure into an underlying tree whose branches are connected by cut bars is a standard procedure. Usually though, the six components of the stress resultant at the cut face are represented by a pair of vectors in 3D, one for the force and one for the moment caused by the stresses in the cut face. A graphical representation of the forces - the so-called Vector-based Graphics Statics~\cite{MaxwellReader} can be obtained by conjoining the various force vectors that meet at each node. This is equilibrium expressed by the closure of the polygon of forces, and in the 3D case, the vectors which form the edges of the polygon do not all need to lie on the same plane. A related construction is also possible using the moment vectors, since nodal equilibrium requires all end moments on all bars meeting at a joint to sum to zero, although that approach has received far less attention.

The key differences hinge on the notions of ``parallel'' and ``perpendicular''. In traditional 2D graphic statics, bar forces can be drawn parallel to their corresponding bars (the so-called Cremona convention) or perpendicular to the bars (the so-called Maxwell convention). Ultimately one obtains the same diagram, but one is a 90 degree rotation of the other. In three dimensions, however, one obtains completely different objects. The ``parallel'' convention corresponds to Vector-based Graphic Statics, with forces denoted as vectors. The ``perpendicular'' convention corresponds to the Rankine construction, with force represented by the area of a flat polygon perpendicular to the bar in question. This paper presents a generalisation of the latter perspective, with assemblages of loops taking over the role of the polyhedra and polytopes of the Rankine approach.

\begin{figure*}[ht!]
\begin{center}
\includegraphics[width=0.8\textwidth]{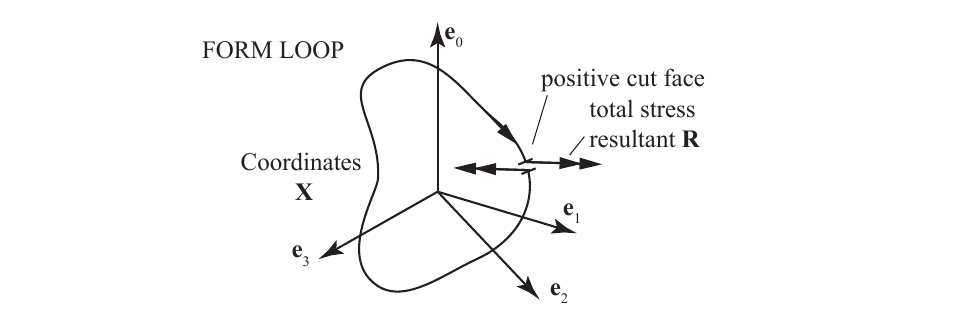}
\end{center}
\caption{Sign convention: each loop has a cut at a point in a non-tree bar, and the orientation of the loop is defined by arbitrarily selecting one of the faces at the cut as the positive face. }
\label{SignConvention}
\end{figure*}

\subsection{Sign Convention}
Each loop is assigned - arbitrarily - an orientation. When a structural bar is (imagined to be) cut, it reveals two cut faces. The two opposing faces correspond to the same loop but with opposing orientations. One face at each cut is selected - arbitrarily - as the positive face, and that then defines which orientation has been selected as positive for the loop that passes through that cut. The positive orientation of the loop is that where the loop's circulation leaves the positive cut face (see Fig.~\ref{SignConvention}). Cuts are made on bars not in the spanning tree, thus only one loop passes through each cut.

\section{Examples}
\subsection{The frame with K5 graph}
Maxwell's classic 1864 paper~\cite{maxwell1864a} begins with the simplest example of an axially self-stressed structure. It has a K4 graph (four nodes, each connected to every other node) and the reciprocal has the same K4 topology. The structure consists of an outer triangle (K3), with three radial spokes connecting to a central node. The K4 is thus the coning of a K3. Another familiar example having this topology is the cross-braced quad - a four-bar linkage with two diagonals which do not connect to each other, and the Maxwell reciprocal diagram is also a cross-braced quad. Each diagram is the 2D projection of a 3D tetrahedron. 

\begin{figure*}
\begin{center}
\includegraphics[width=0.8\textwidth]{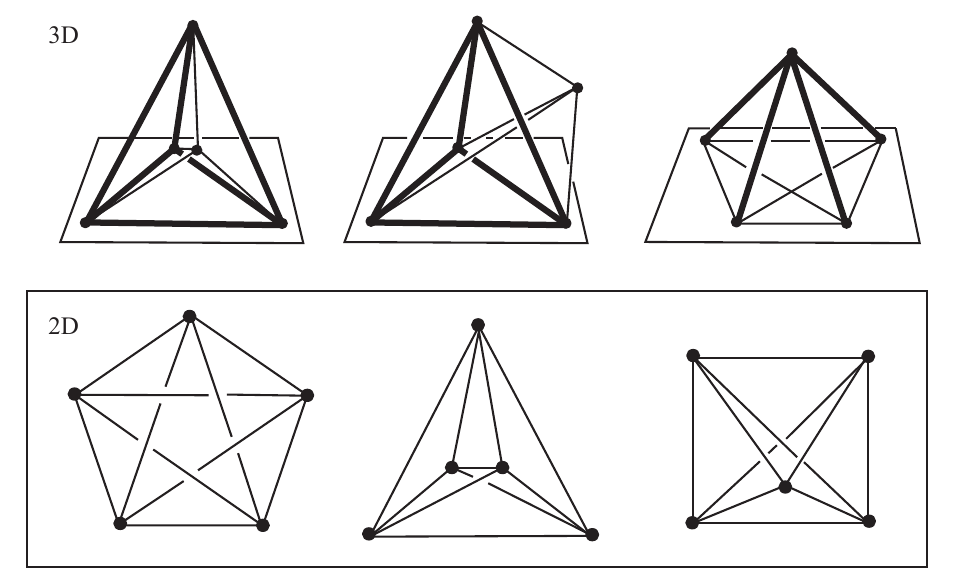}
\end{center}
\caption{A variety of frames with the K5 graph. Some are fully 3D and some fully 2D. Top right is the 3D coning of a 2D K4.
}
\label{K5varieties}
\end{figure*}

Moving everything up a dimension, the simplest example of an axially self-stressed structure in 3D has the K5 topology: there are five nodes, each connected to every other (see Fig.~\ref{K5varieties}). When the structure is a 3D truss in general position, the Rankine reciprocal (where forces are represented by area) has the same K5 topology.
The structure consists of an outer tetrahedron (K4) with four radial spokes connecting to a central hub. The K5 is thus the coning of a K4. Again, other geometries are possible which have the same topology. For example,  the central node may be moved outside the original outer tetrahedron. More subtly, the central node may be moved to lie in one of the original outer faces of the tetrahedron. This gives a plane K4 substructure, which is coned to a single out-of-plane node (see Fig.~\ref{K5varieties}, top right). Such 3D configurations which involve 2D substructures are commonplace in engineering, but have previously been difficult to deal with in graphic statics. Similarly, all five nodes may lie in the same plane. Care with terminology is needed: the five nodes lie in a plane but the K5 graph is not a {\it planar} graph - it cannot be drawn on a plane without a crossing. 
One consequence of this is that the 2D K5 structure has no 2D Maxwell reciprocal figure. However, a reciprocal figure of sorts can be obtained via the usual work-around of adding a node at any crossing in the form diagram.

Each of the K5 figures in Fig.~\ref{K5varieties} is the projection (to 3D or 2D) of a 4D polytope, a 5-cell. A 5-cell consists of five tetrahedra, each pair of which shares a triangular face. At this stage of the development of the theory, however, the concept of 3D volumetric cells such as tetrahedra is not required. Here, loops are the highest-order objects. 

 The dual K5 figures provide the simplest example of 
 Maxwell-Rankine 3D reciprocals. Each triangular face in the reciprocal is normal to a corresponding bar in the form diagram, with the area of the triangle giving the magnitude of the bar force. Before illustrating more general reciprocals for a K5 frame (which may involve moments, for example) we first show how the traditional Maxwell-Rankine reciprocal emerges naturally from the new homology-based description. 
\begin{figure*}[ht!]
\begin{center}
\includegraphics[width=\textwidth]{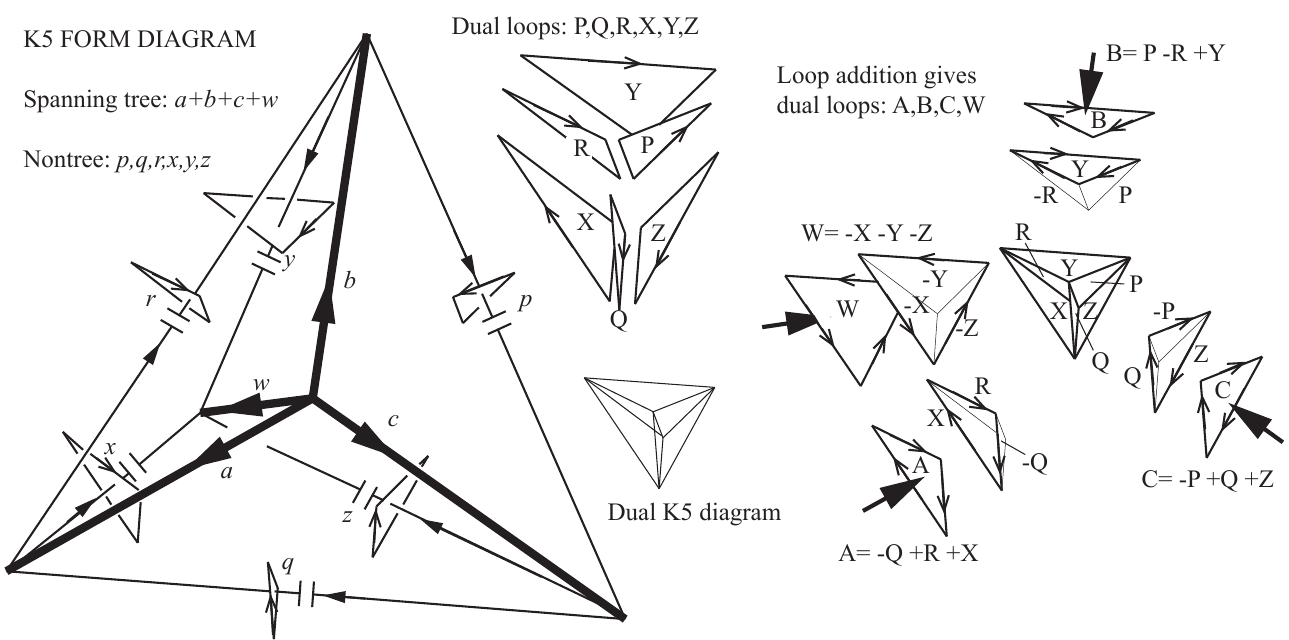}
\end{center}
\caption{a) The frame with K5 graph has ten bars, four of which can form a spanning tree (bars $a,b,c,w$ here).  Each of the six non-tree bars $(p,q,r,x,y,z)$ defines a loop. A purely axial solution may be represented by six triangular loops $(P,Q,R,X,Y,Z)$, one dual to each form loop. b) These six dual triangles may be arranged to share a central node. c) Loops representing forces in tree bars (such as bar $a$) may be determined by summing the three dual loops ($A = -Q + R+ X$ here) appropriate to that bar. The summation uses the group operation in the edge group $C_1$.}
\label{K5}
\end{figure*} 

\subsubsection{The standard axial solution for the K5 truss revisited}
A 3D frame with K5 topology is shown in Fig.~\ref{K5}a.  
The first step in the new description assigns an orientation to each bar. These are indicated by arrows in the diagram. For each bar, any choice of orientation is admissible.
Next, a spanning tree is selected. Many spanning trees are possible, and any would be admissible. Here we choose the bars $a,b,c,w$ that radiate from the central node to the four corners of the surrounding tetrahedron. 
Each of the six non-tree bars $p,q,r,x,y,z$ defines a loop that connects through the tree. For example, in Fig.~\ref{K5}, the non-tree bar $p$ defines the loop $p - c + b$.

The reciprocal figure thus consists of six dual loops $P,Q,R,X,Y,Z$ (see Fig.~\ref{K5}b). Since there are no moments in the axial solution, each reciprocal loop can be a simple polygon in 3D normal to its corresponding bar, and for simplicity, we can let each polygon be a triangle. Each tree bar belongs to three loops, and we require each such force to be axial. This can be readily achieved if the non-tree loops are arranged so that the sum of any three adjacent triangles has a triangular boundary normal to a tree bar (see Fig.~\ref{K5}c). Loops are said to be ``adjacent'' here if they share a tree bar. For example, loops defined by bars $p$, $y$ and $r$ are adjacent as tree bar $b$ is common to all three.

Each tree bar ($b$, say) is an element of three structural loops ($p,y,r$ here).
The force in each such tree bar is then given by the sum of the three corresponding dual force loops ($P,Y,R$ here) (see Fig.~\ref{K5}c).  These three force loops may be chosen to be triangles whose bivector sum is a triangle $B$ normal to $b$. Addition is all performed with respect to the homological algebra. Thus, in this example, the positive face of a cut on bar $b$ must carry a force $B = P+Y-R$, all chosen such that $B$ is normal to $b$. 
%There are four sets (one per tree bar) of three adjacent loops, and the duals of each set of loops can be arranged as three faces of a tetrahedron. The missing fourth face will be dual to the tree bar. For example, tree bar $a$ is a member of loops defined by $q$, $r$ and $x$. With the orientations of bars and loops as illustrated in Fig.~\ref{K5}, the loop $A$ dual to $a$ 
%is given by $A = -Q+R+X$, with the operations of addition and subtraction of loops as defined in the homological algebra. The result, illustrated schematically in Fig.~\ref{K5}c, is the oriented loop $A$ whose bivector area corresponds to the force in bar $a$. 
If desired, any force bivector in 3D can be converted to a force vector in 3D by taking the Hodge dual, 
i.e.~the force in bar $b$ may be represented by a vector normal to the triangle $B$, of magnitude equal to the triangle area and with the force direction on the positive cut face of bar $b$ given by the orientation of the dual loop $B$, taking a right-hand screw rule. 

Collecting all the triangular force loops in this manner, we arrive at a geometric object which has essentially the same geometry as the traditional Maxwell-Rankine solution for a 3D K5 truss. However, rather than being a set of 5 adjacent tetrahedra with 10 triangular faces (i.e.~a 5-cell), the reciprocal object consists of only the six triangular loops $P,Q,R,X,Y,Z$  which all share a corner at the centre of the diagram. The outer four triangles $A,B,C,W$ are defined implicitly, as described above. Moreover, there is as yet no concept of ``face'' or ``cell'' or ``polytope''. The next example shows how the state of axial self-stress in a simple tensegrity may be represented within the loop formalism.

%in that their actions are represented by the bivector sum of three adjacent triangles of the basic six. Fig.~\ref{K5}c shows the construction in more detail, paying attention to the sign convention. The orientation of a dual loop (using a right-hand screw rule along the outward normal to the cut face of the bar) then defines whether it denotes tension or compression.

\subsection{Axial forces in simple tensegrities}
Simple tensegrity structures possess states of pure axial self-stress which can
be difficult to represent with traditional Maxwell-Rankine reciprocals. In many cases, the tensegrity has the topology of a spherical polyhedron, and for the usual 3D form-force duality (where cells are dual to points), the reciprocal object would thus appear to consist of only a single point (or perhaps more strictly, two points, if the form topology is considered to be the polytope which is the double cover of the structural polyhedron). A reciprocal object with only one or two points does not have sufficient richness to represent the many internal bar forces. Another problem is that some tensegrities (such as the one with cube topology considered later) have polygonal `faces'  which are not plane, and these are difficult to handle in most other approaches to graphic statics.

Such problems have been considered at length previously (see, for example,~\cite{ McRobieHamburg, McRobieBoston}). One solution involves coning the polyhedron, with all structural nodes being connected to an inner pole via so-called Zero Bars~\cite{McRobieZero}. (A Zero Bar is one whose reciprocal force polygon has zero oriented area). The coning creates a set of internal pyramidal cells, each with apex at the pole of the coning, with the base of each pyramid coinciding with a face of the surrounding polyhedron.  An alternative construction augments each face of the polyhedron with a ``face cushion'' cell~\cite{McRobieHamburg}. A face cushion is a cell with only two faces, and these share a single surrounding boundary edge loop.  A stress function may be associated with each face cushion is such a way as to define extra reciprocal nodes, and the greater number of reciprocal nodes allows the creation of a reciprocal object of sufficient richness to represent the internal forces. This construction is described in a later section. Here, though,  we give an example of a simple tensegrity whose state of self-stress is described by the methods of this paper, using only loops. 

\subsubsection{The octahedral three-prism tensegrity}

\begin{figure*}[b]
\begin{center}
\includegraphics[width=\textwidth]{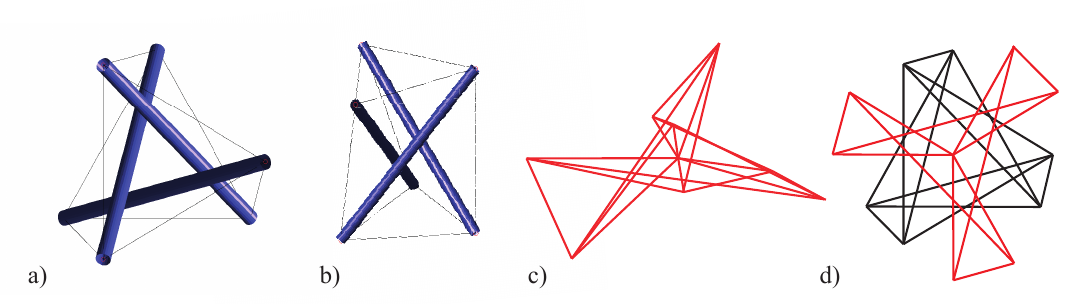}
\end{center}
\caption{The three-prism tensegrity with octahedral topology. a) Plan view. b) General view. c)~Rankine reciprocal. d) Plan view of structure and reciprocal.  
}
\label{PrismFirst}
\end{figure*} 

\begin{figure*}[ht]
\begin{center}
\includegraphics[width=\textwidth]{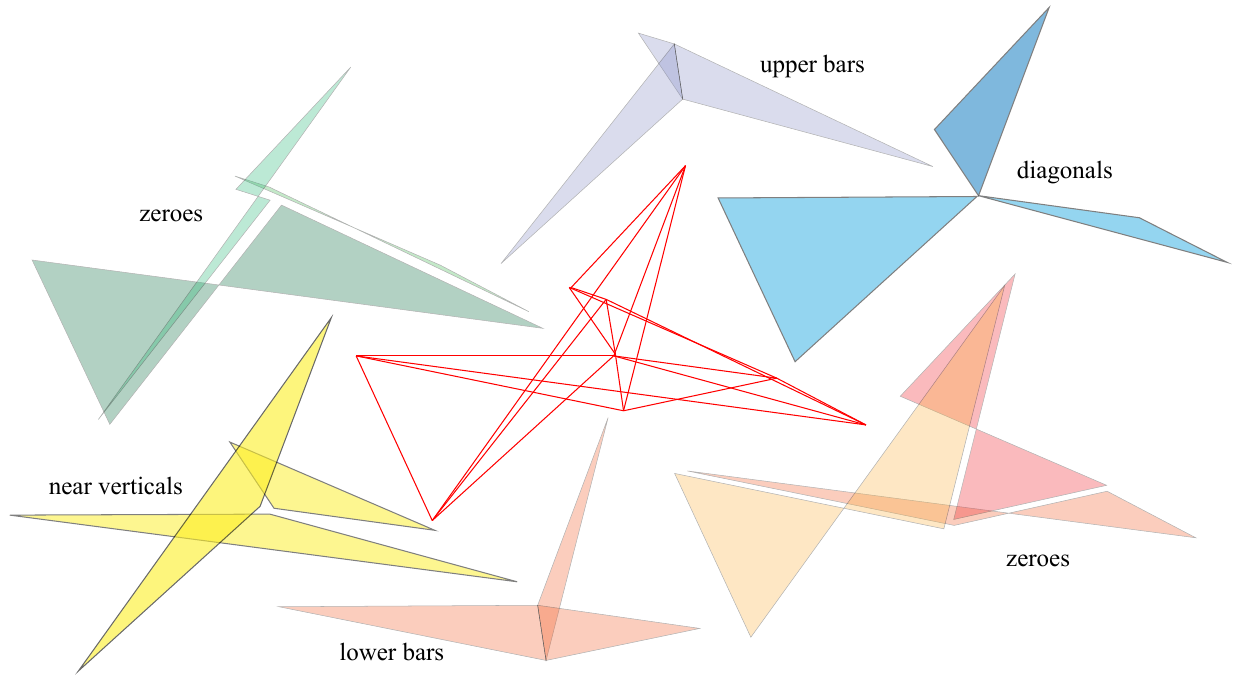}
\end{center}
\caption{The decomposition of the Rankine reciprocal for the three-prism tensegrity into the individual polygons dual to each structural bar of Fig.~\ref{PrismFirst}b). The forces in the main compression members are represented by the blue triangles (top right). The forces in the upper and lower triangles of tension members are shown in purple (upper) and pink (lower). The forces in the near-vertical tension members are given by the yellow triangles. The three green (left) and three orange (right) quads each have zero area, corresponding to lack of force in the zero bars that cone the structure to a central node. 
}
\label{PrismDecomposition}
\end{figure*}

The three-prism tensegrity with octahedral topology is shown in Fig.~\ref{PrismFirst}, together with the Rankine reciprocal that results from a Zero Bar coning to a central node, as per McRobie~\cite{McRobieHamburg}. The reciprocal is decomposed into its constituent polygons in Fig.~\ref{PrismDecomposition}. % The Rankine-Minkowksi diagram, where structural members are augmented with their reciprocal polygons, is shown in Fig.~\ref{PrismMink}

\begin{figure*}[hp]
\begin{center}
\includegraphics[width=\textwidth]{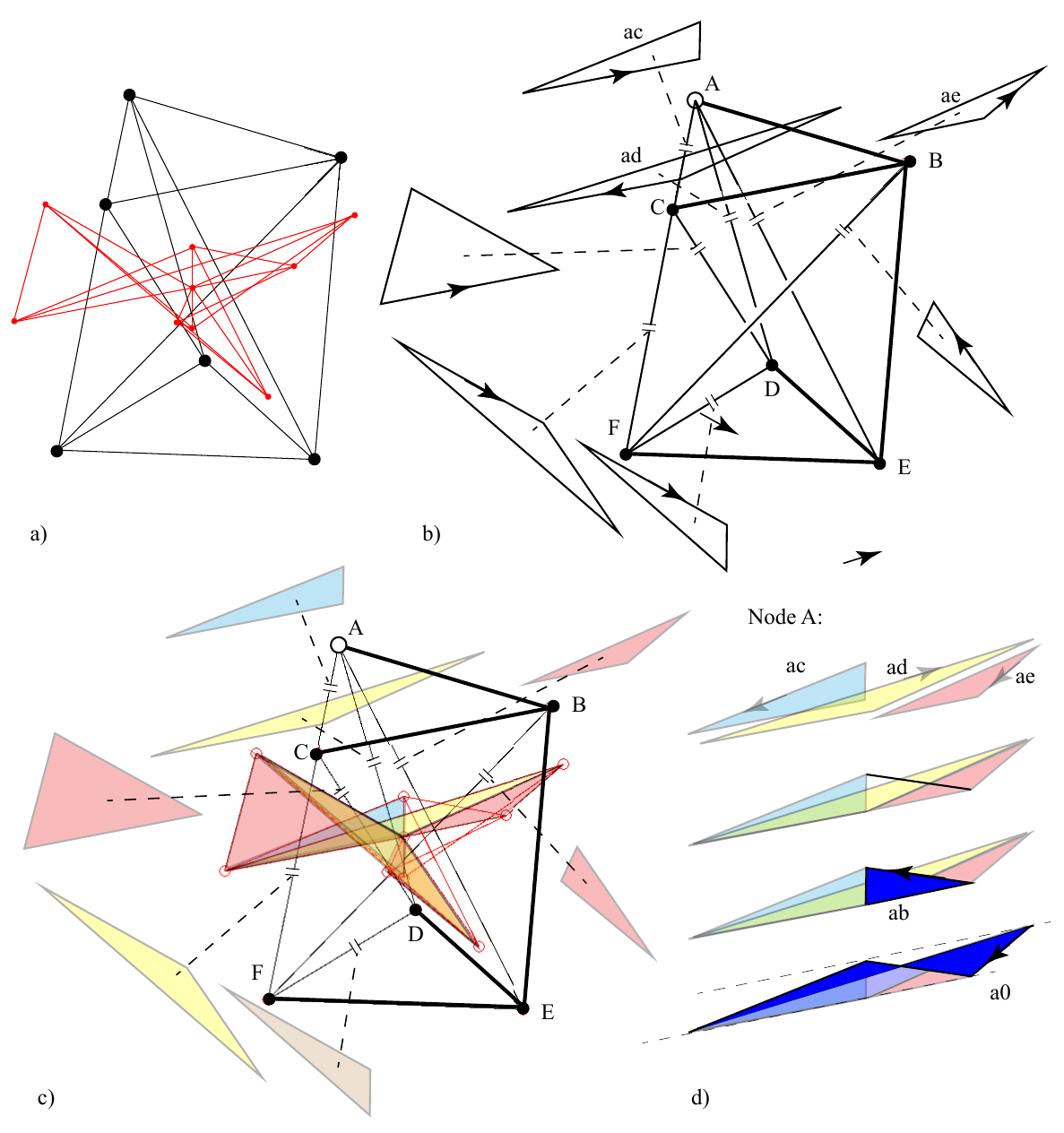}
\end{center}
\caption{a) The three-prism tensegrity (black) and a Rankine reciprocal (red). b) A spanning tree with seven non-tree bars, to each of which is ascribed a simple triangular dual loop. and seven loops These six dual triangles may be arranged to share a central node. c) The seven reciprocal loops are faces of the reciprocal polyhedron. d) Forces in tree bars are represented by appropriate sums of oriented loops. Equilibrium requires a closed polyhedron at Node A, which is a pyramid having Zero Face $\mathbf{a0}$ as base, and triangular faces $\mathbf{ab}$, $\mathbf{ad}$, $\mathbf{ae}$ and $\mathbf{ab}$.
}
\label{Prism1}
\end{figure*}

Here we illustrate how such objects fit within the loop-to-loop description. Fig.~\ref{Prism1}a shows the so-called three-prism tensegrity (in black), together with its Rankine reciprocal (in red). The structure has the topology of an octahedron, and as described in~\cite{McRobieHamburg}, a Rankine reciprocal may be obtained via coning to a central node connected by Zero Bars which triangulate the interior of the surrounding octahedron. 

Fig.~\ref{Prism1}b shows the loop-to-loop description, which does not require the addition of Zero Bars. The structure has six nodes, and one possible spanning tree of 5 bars is shown, consisting of the bars BE, BA, BC, ED and EF. 

There are seven non-tree bars, each of which defines a structural loop. Each non-tree bar may be cut, and self-equilibrating forces and moments may be applied across each cut. If the structure were a fully-welded frame, there would thus be $7\times 6 = 42$ possible states of self-stress (six per cut). Given that each node has valence four and the bars connecting at a node do not contain any three coplanar bars, it follows from elementary structural analysis that there is (at most) only one state of self-stress which is purely axial. (If we know one force at any node, then the other three are uniquely determined by equilibrium. We may then move to an adjoining node, and determine the unknown forces there. If the procedure closes, then there is a state of self-stress. If it does not, then there is no state of self-stress). 

Dual to each structural loop a reciprocal force loop may be assigned in order to define the forces and moments at the cuts. Here, these reciprocal loops have been chosen to be simple triangles (Fig.~\ref{Prism1}b). Since there are no bending or torsional moments involved in the purely axial solution,  we focus on how the reciprocal triangles are located in 3D space. Within 3D, each reciprocal triangle must be perpendicular to the bar whose force it represents.  

Here, we shall use the results of the preceding analysis that used Zero Bars (Fig.~\ref{PrismFirst}c). We are thus limiting the presentation to showing how the state of axial self-stress may be represented with loops, rather than tackling the more general problem which would require finding an axial-only reciprocal given only the structural geometry.

(Aside: a nice geometry for this configuration locates the six nodes on a cylinder of unit radius. The cylinder is oriented along the $z$ axis, with top and bottom planes at $z = \pm 0.5$. The Maxwell-Rankine stress function is set to $f = 0.25$ at all nodes, and the structure is coned to a pole at the origin where $f=0$. The coning creates a polytope consisting of 8 tetrahedra within an overall octahedral hull. 
Via the Legendre transform, the stress function gradients define the dual object which has 9 nodes, six of which lie on the unit cylinder and three on the $z$-axis. The dual stress function takes the value $\phi= 0$ at all nodes, except the node at the origin, where $\phi = -1$.)

\begin{figure*}[h!]
\begin{center}
\includegraphics[width=\textwidth]{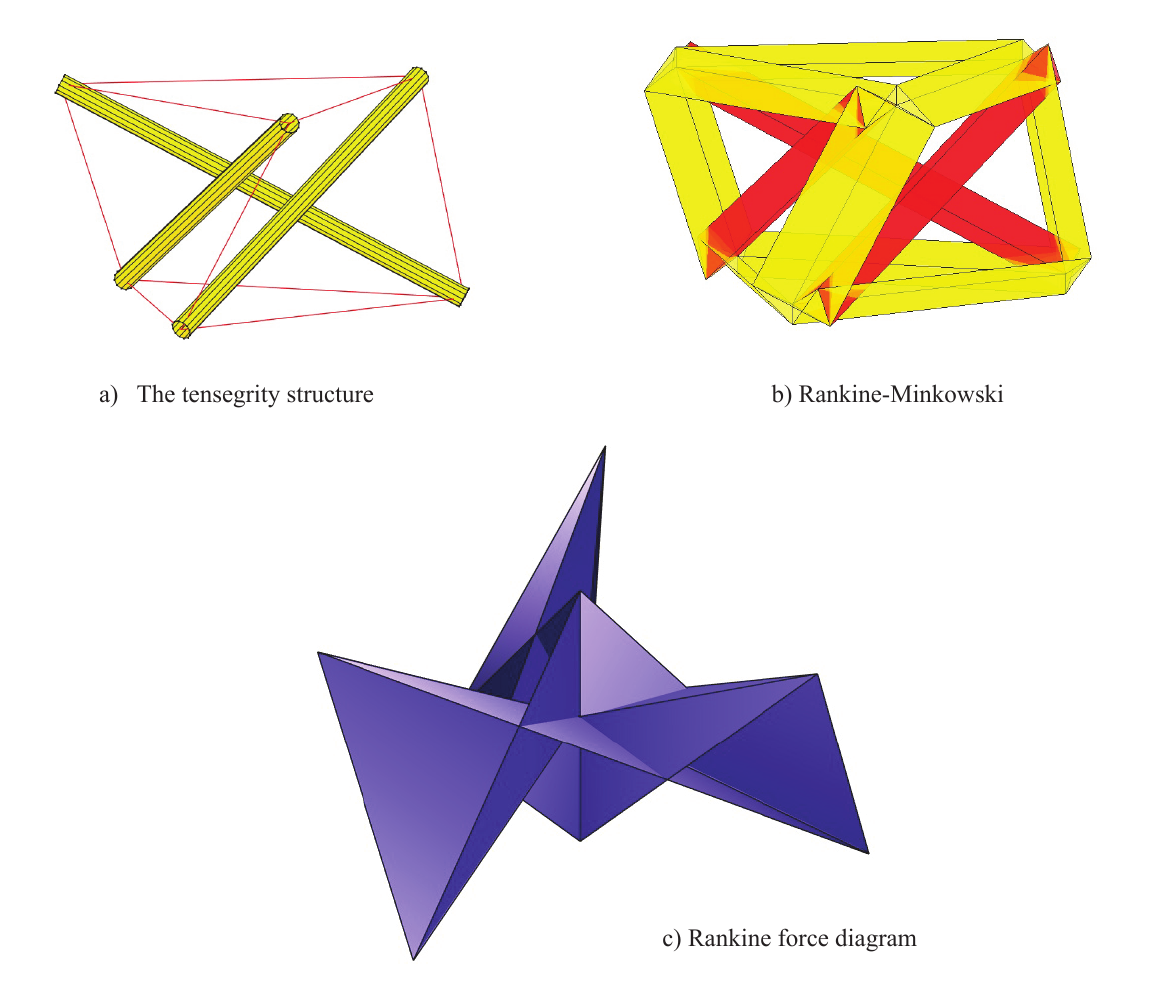}
\end{center}
\caption{The three-prism tensegrity with octahedral topology. a) Structure. b) Rankine-Minkowski. c) Rankine force diagram for the state of axial self-stress.
}
\label{Combined}
\end{figure*}

Returning to the loop formalism, the force in any tree bar may be obtained by adding all oriented force loops whose form loops involve that tree bar. For example,   Fig.~\ref{Prism1}d shows how the force in bar $ab$ may be obtained. Equilibrium at node A requires the forces in all four bars there to be in equilibrium. Equilibrium may be represented by a closed polyhedron, and in this case, rather than being a simple tetrahedron, we have chosen a pyramid with five faces, consisting of four triangular sides which surround a zero-area base $\mathbf{a0}$. This zero-area base lies in a plane normal to the line connecting the corner A to the point at the centre of the structure. This zero-area base is a ghost-like echo of the Zero Bar in the earlier analysis which connected corner A to the fictitious central node. Note however, that no Zero Bars are required in the current analysis. To be explicit, the oriented areas of the triangular force polygons of the bars meeting at A add to zero, corresponding to equilibrium, but they do not form a closed polyhedron. To close the polyhedron, a fifth polygon would be required, but this polygon has zero oriented area.

(Aside: by focusing on the projections of the reciprocal loops in the 3D space of $(x,y,z)$ we have ignored their projections on the planes which involve the stress function dimension. Despite this omission, a sensible set of axial forces was obtained. However, the loops do not lie completely in $(x,y,z)$. If they did, this would mean the total moment on any cut face would be zero. However, the total moment is the sum of the internal bending and torsional moments together with the moment created by the bar force about the origin. Since bars do not pass through the origin, then loops in 3D alone would require a set of bending moments and torsions to equilibrate with the bar-force-induced moments. It follows that the reciprocal object corresponding to the axial-only solution cannot live purely in 3D $(x,y,z)$. As described in the earlier aside, one possibility is that the dual stress function $\phi$ is zero at all reciprocal nodes except the origin.)

Many alternative constructions are admissible. At each of the cuts on the seven non-tree bars, any polygon of the correct oriented area can assigned, with ``correct'' meaning that the oriented areas are perpendicular to all twelve bars, and that oriented areas sum to zero at each node.

%One example would be to choose circles normal to each cut bar, with areas matching the bar force magnitudes. The forces in tree bars would then be given by the bivector sum of three or four oriented circles. 
%The resultant could then be represented by an oriented circle. Extruding the circles along each bar would create a graphic indicative of required bar areas, not unlike the post-processing graphics of many finite element packages. However, this loses much of the elegance of the Maxwell-Rankine construction, where equilibrium is represented by closed polyhedra ``by hydrostatics'', to use Maxwell's phrase~\cite{maxwell1870a}. Indeed,  the elegance of Fig.~\ref{PrismFirst}d suggests that there is something natural about the choice of the reciprocal polygons that result from the Zero Bar coning construction. This is even more pronounced in Fig.~\ref{PrismFirst}c.

The results of the graphical analysis of the three-prism tensegrity are summarised in Fig.~\ref{Combined}, showing the structure, the Rankine-like reciprocal for the state of axial self-stress, and the two diagrams combined into the Rankine-Minkowski diagram.  (It is ``Rankine-like'' as, strictly speaking, there is no Rankine traditional Rankine reciprocal, and the bow-tie-shaped "faces" of zero area are only implicitly defined within the assemblage of triangular loops). 
The material point is that this a highly graphical approach.

\begin{figure*}[b!]
\begin{center}
\includegraphics[width=\textwidth]{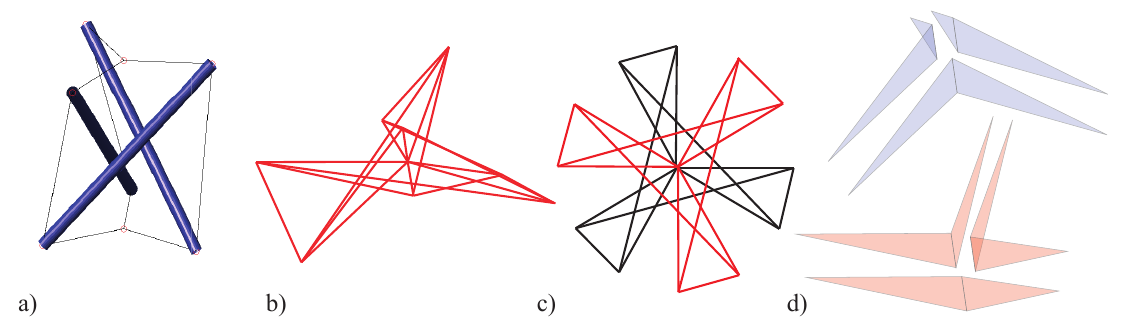}
\end{center}
\caption{The three-prism tensegrity with cubic topology. a) General view. b) Rankine reciprocal as previously. c) Plan view of structure and reciprocal. d) Objects reciprocal to the radial star members.
}
\label{PrismSecond}
\end{figure*} 

Closely related to the octahedral three-prism tensegrity is the cube tensegrity, where the upper and lower triangles have been replaced with the ``Mercedes" star configuration, each with a central node on the $z$-axis. There are some interesting dualities here. Obviously, a cube is the Platonic dual of an octahedron. More interestingly the reciprocal of an octahedral three-prism tensegrity in Vector-based Graphic Statics is a cube tensegrity. In the Rankine-like loop formalism here, though, the cube and octahedral tensegrities are not dual to each other, Instead, each is dual to a more complicated object (Fig.~\ref{PrismDecomposition}). These force diagrams share a set of nodes and lines, but the decomposition into loops differs in the two cases. Dual to each edge of the upper and lower triangular faces of the octahedron there is a triangular loop  (Fig.~\ref{PrismDecomposition}). However, dual to each spoke of the Mercedes symbol in the cube configuration, there are two such triangular loops whose sum within the homological algebra is a bent quad (Fig.~\ref{PrismSecond}d).

\section{Summary}
This paper has shown how any frame structure can decomposed into a set of loops and any state of self-stress within the structure can be represented by a set of dual loops in the 4D force space. The projected areas of the dual loops give the six components of the forces and total moments about the origin at any cut face on a structural loop. Where Part 1 of this sequence of papers focused on how moments could be represented, with simple examples that typically involved only a single loop,
the paper here focused instead on how frames involving multiple loops may be modelled. 

Although all the apparatus for modelling moments is included in the formalism, this paper restricted its attention - for simplicity of explanation - to examples having purely axial states of self-stress.  The examples, a structure with a K5 graph and a prism tensegrity, constitute the archetypal first steps in any approach to 3D graphic statics, and it was shown how each could be treated with the loop formalism. The tensegrity example is somewhat non-trivial in the usual 3D Rankine approach, given that it does not possess a standard Rankine reciprocal.  

The highest dimensional objects used in the formalism presented here are loops. It has not yet been suggested that the loops might be the boundaries of polygonal faces or more general surfaces, or that the polygonal faces might be the boundaries of cells. (Strictly speaking, though, the notions of surfaces and cells were used offline to find the Rankine reciprocals that were then subsequently presented in terms of loops. The loop approach leads to infinitely many possible reciprocals, but for the chosen examples, the Rankine reciprocals have a pleasingly neat geometry.)  In the papers that follow, the homological algebra is extended to include surfaces, cells and higher dimensional objects.
The key elements will be CW-complexes, which provide a convenient generalisation of polygons, polyhedra and polytopes. The formalism will also be extended to describe kinematics.
However, the main points of the present paper are that, despite the simplicity of the examples presented, the approach can deal with the forces and moments in any 3D frame, and that  - as exemplified by the figures throughout - the approach is unapologetically graphical.

\bibliography{McRobieBib2}
%\printbibliography
\end{document}